\documentclass[12pt]{amsart}
\usepackage{amsfonts}
\usepackage{ifthen}
\usepackage{amsthm}
\usepackage{amsmath}
\usepackage{graphicx}
\usepackage{amscd,amssymb,amsthm}
\usepackage{graphicx}
\usepackage{epstopdf}
\usepackage{hyperref}
\usepackage{color}

\newcounter{minutes}
\divide\time by 60
\newcounter{hours}
\multiply\time by 60 \addtocounter{minutes}{-\time}

\newtheorem{lemma}{Lemma}[section]
\newtheorem{theorem}{Theorem}[section]


\title{Radius of Uniform Convexity of some special functions}

\author[\.{I}. Akta\c{s}]{\.{I}brah\.{I}m Akta\c{s}}
\address{Department of Mathematical Engineering, Faculty of Engineering and Natural Sciences, G\"{u}m\"{u}\c{s}hane University, G\"{u}m\"{u}\c{s}hane, Turkey}
\email{aktasibrahim38@gmail.com}

\author[E. Toklu]{Evr{\.I}m Toklu}
\address{Department of Mathematics, Education Faculty, A\u{g}r{\i} {\.I}brah{\.I}m \c{C}e\c{c}en University, 04100 A\u{g}r{\i}, Turkey} \email{etoklu@agri.edu.tr}

\author[H. Orhan]{Hal\.{i}t Orhan}
\address{Department of Mathematics, Faculty of Science, Atat\"{u}rk University, Erzurum, Turkey}
\email{orhanhalit607@gmail.com}

\begin{document}

\maketitle

\begin{abstract}
In this investigation our main aim is to determine the radius of uniform convexity of the some normalized $ q $-Bessel and Wright functions. Here we consider six different normalized forms of $ q $-Bessel functions, while we apply three different kinds of normalizations of Wright function. Also, we have shown that the obtained radii are the smallest positive roots of some functional equations.
\end{abstract}

\section{Introduction and Preliminaries}
Special and geometric function theories are the most important branches of the mathematical analysis. There is a close relationship between special and geometric functions theory since hyper geometric functions has been used in the proof of famous Bieberbach conjecture. This is why most of mathematicians have considered the some geoemetric properties of special functions which can be expressed by the hyper geometric series. Especially, the some geometric properties of Bessel, Struve, Lommel, Wright and $ q $-Bessel functions have been investigated by the many authors. The first important results on the geometric properties of hyper geometric and related functions can be found in \cite{brown,todd,Merkes,wilf}. Actually, there are some relationships between the geometric properties and the zeros of the special functions. Owing to these relations numerous investigations has been done on the zeros of the above mentioned special functions. The comprehensive informations about Bessel function and its $q$-analogue can be found in \cite{Wat}, and also the some results on the zeros of some special functions can be found in \cite{bsingh,ismail1,ismail2,koelink,Koornwinder,steinig1,steinig2}. Recently, some geometric properties (like univalence, starlikeness, convexity and uniform convexity) of Bessel, Struve and Lommel functions of the first kind have been investigated by \cite{aktas2,aktas1,bdoy,bks,btk,bos,samy,basz,BY,Deniz,szasz,wright}. In addition, the radii of starlikenes and convexity of some normalized $q$-Bessel functions have been studied in \cite{aktas3,aktas4,BDM}. Motivated by the previous works in this field our aim is to determine the radius of uniform convexity of the some normalized $q$-Bessel and Wright functions.

Now, we would like to present the some basic concepts related with geometric function theory. 
Let $\mathbb{D}_r$ be the open disk $\{z\in\mathbb{C}:\left|z\right|<r\}$ with radius $r>0$ and $\mathbb{D}_1=\mathbb{D}$. Let $\mathcal{A}$ denote the class of analytic functions $f:\mathbb{D}_r\rightarrow\mathbb{C},$
\begin{equation}\label{1.1}
f(z)=z+\sum_{n\geq2}a_{n}z^n,
\end{equation}
which satisfy the normalization conditions $f(0)=f^{\prime}(0)-1=0$. By $\mathcal{S}$ we mean the class of functions belonging to $\mathcal{A}$ which are univalent in $\mathbb{D}_r$.  On the other hand, the class of convex functions is defined by $$\mathcal{K}=\bigg\{f\in \mathcal{S}:\Re\left(1+\frac{zf^{\prime \prime}(z)}{f^{\prime}(z)}\right)>0 \text{ for all } z\in \mathbb{D}_r \bigg\}.$$ The radius of convexity of an analytic locally univalent function $f:\mathbb{C}\rightarrow\mathbb{C}$ is defined by $$r^{c}(f)=sup\bigg\{r>0:\Re\left(1+\frac{zf^{\prime \prime}(z)}{f^{\prime}(z)}\right)>0 \text{ for all } z\in \mathbb{D}_r \bigg\}.$$ Note that $r^{c}(f)$ is the largest radius for which the image domain $f\left(\mathbb{D}_{r^{c}(f)}\right)$ is a convex domain in $\mathbb{C}.$ For more information about convex functions we refer to Duren's book \cite{Duren} and to the references therein.

In \cite{Goodman} the author has introduced the concept of uniform convexity for the functions of the form \eqref{1.1}.	A function $f(z)$ is said to be uniformly convex in $\mathbb{D}$ if $f(z)$ is in class of usual convex functions and has the property that for every circular arc $\gamma$ contained in $\mathbb{D}$, with the center $\zeta$ also in $\mathbb{D}$, the arc $f(\gamma)$ is a convex arc.
An analytic description of the uniformly convex functions has been given by R\o{}nning in \cite{Ronning} read as follows:
\begin{theorem}
	Let $f(z)$ is of the form \eqref{1.1}. Then, $ f $  is a uniformly convex functions if and only if 
	\begin{equation*}
	\Re\left(1+\frac{zf^{\prime\prime}(z)}{f^{\prime}(z)}\right)>\left|\frac{zf^{\prime\prime}(z)}{f^{\prime}(z)}\right|,    z\in\mathbb{D}.
	\end{equation*}
\end{theorem}
On the other hand, the concept of the radius of uniform convexity is defined by (see \cite{Deniz}) 
$$r^{uc}(f)=\sup\Big\{r\in(0,r_f): \Re\left(1+\frac{zf^{\prime\prime}(z)}{f^{\prime}(z)}\right)>\left|\frac{zf^{\prime\prime}(z)}{f^{\prime}(z)}\right|, z\in\mathbb{D}\Big\}.$$
Thanks to above theorem we can determine the radius of uniform convexity for the functions of the form \eqref{1.1}. Also, we will need the following lemma in the sequel.

\begin{lemma}[\cite{Deniz}]
If  $a>b>r\geq \left|z\right|,$ and $ \lambda\in [0,1] $, then
\begin{equation}\label{1.2}
\left|\frac{z}{b-z}-\lambda\frac{z}{a-z}\right|\leq\frac{r}{b-r}-\lambda\frac{r}{a-r}.
\end{equation}
The followings are very simple consequences of this inequality
\begin{equation}\label{1.3}
\Re\left(\frac{z}{b-z}-\lambda\frac{z}{a-z}\right)\leq\frac{r}{b-r}-\lambda\frac{r}{a-r}
\end{equation}
and 
\begin{equation}\label{1.4}
\Re\left(\frac{z}{b-z}\right)\leq\left|\frac{z}{b-z}\right|\leq\frac{r}{b-r}.
\end{equation}
\end{lemma}

\section{Radius of Uniform Convexity of Some Special Functions}
\setcounter{equation}{0}
In this section we focus on some normalized $q$-Bessel and Wright functions and determine the radii of uniform convexity of these functions.

\subsection{Uniform Convexity of some normalized q-Bessel Functions}
Jackson's second and third (or Hahn-Exton) $q$-Bessel functions are defined as follow:
$$J_{\nu}^{(2)}(z;q)=\frac{(q^{\nu+1};q)_{\infty}}{(q;q)_{\infty}}\sum_{n\geq0}\frac{(-1)^{n}\left(\frac{z}{2}\right)^{2n+\nu}}{(q;q)_{n}(q^{\nu+1};q)_{n}}q^{n(n+\nu)}$$ and $$J_{\nu}^{(3)}(z;q)=\frac{(q^{\nu+1};q)_{\infty}}{(q;q)_{\infty}}\sum_{n\geq0}\frac{(-1)^{n}z^{2n+\nu}}{(q;q)_{n}(q^{\nu+1};q)_{n}}q^{\frac{1}{2}{n(n+1)}},$$ where $z\in\mathbb{C},\nu>-1,q\in(0,1)$ and $$(a;q)_0=1,	(a;q)_n=\prod_{k=1}^{n}\left(1-aq^{k-1}\right),	 (a,q)_{\infty}=\prod_{k\geq1}\left(1-aq^{k-1}\right).$$
It is known that the Jackson's second and third $q$-Bessel functions are $q$-extensions of the classical Bessel function of the first kind $J_{\nu}$. Clearly, for fixed $z$ we have $J_{\nu}^{(2)}\left((1-z)q;q\right)\rightarrow{J}_{\nu}(z)$ and $J_{\nu}^{(3)}\left((1-z)q;q\right)\rightarrow{J}_{\nu}(2z)$ as $q\nearrow1.$

Because the functions $J_{\nu}^{(2)}(.;q)$ and $J_{\nu}^{(3)}(.;q)$ do not belong to $\mathcal{A}$, first we consider the following six normalized forms as in \cite{BDM}. For $\nu>-1$,
$$f_{\nu}^{(2)}(z;q)=\left(2^{\nu}c_{\nu}(q)J_{\nu}^{(2)}(z;q)\right)^{\frac{1}{\nu}}, \nu\neq0$$
$$g_{\nu}^{(2)}(z;q)=2^{\nu}c_{\nu}(q)z^{1-\nu}J_{\nu}^{(2)}(z;q),$$
$$h_{\nu}^{(2)}(z;q)=2^{\nu}c_{\nu}(q)z^{1-\frac{\nu}{2}}J_{\nu}^{(2)}(\sqrt{z};q),$$
$$f_{\nu}^{(3)}(z;q)=\left(c_{\nu}(q)J_{\nu}^{(3)}(z;q)\right)^{\frac{1}{\nu}}, \nu\neq0$$
$$g_{\nu}^{(3)}(z;q)=c_{\nu}(q)z^{1-\nu}J_{\nu}^{(3)}(z;q),$$
$$h_{\nu}^{(3)}(z;q)=c_{\nu}(q)z^{1-\frac{\nu}{2}}J_{\nu}^{(3)}(\sqrt{z};q),$$ where $c_{\nu}(q)=(q;q)_{\infty}\big/(q^{\nu+1};q)_{\infty}$. Consequently, all of the above functions belong to the class $\mathcal{A}$. Of course there exist infinitely many other normalization for both Jackson and Hahn-Exton $q$-Bessel functions, the main motivation to consider these six functions is the fact that their limiting cases for Bessel functions appear in literature, see for example \cite{brown} and the references therein.

It is known from \cite[Lemma 1., p.972]{BDM} that, if $\nu>-1$ then the Hadamard factorizations of the functions $z\mapsto J_{\nu}^{(2)}(z;q)$ and $z\mapsto J_{\nu}^{(3)}(z;q)$ are of the form
\begin{equation*}
J_{\nu}^{(2)}(z;q)=\frac{z^{\nu}}{2^{\nu}c_{\nu}(q)}\prod_{n\geq1}\left(1-\frac{z^2}{j_{\nu,n}^2(q)}\right)
\end{equation*}
and
\begin{equation*}
J_{\nu}^{(3)}(z;q)=\frac{z^{\nu}}{c_{\nu}(q)}\prod_{n\geq1}\left(1-\frac{z^2}{l_{\nu,n}^2(q)}\right)
\end{equation*}
where $j_{\nu,n}(q)$ and $l_{\nu,n}(q)$ are the $ n $th positive zeros of the functions $ J_{\nu}^{(2)}(z;q) $ and $J_{\nu}^{(3)}(z;q)$.
Also, it is known from \cite[Lemma 7., p.975]{BDM} that, if $\nu>0$ then the Hadamard factorizations of the derivatives of the functions $z\mapsto J_{\nu}^{(2)}(z;q)$ and $z\mapsto J_{\nu}^{(3)}(z;q)$ are of the form 
\begin{equation*}
\frac{dJ_{\nu}^{(2)}(z;q)}{dz}=\frac{\nu {\left(\frac{z}{2}\right)}^{{\nu}-1}}{2c_{\nu}(q)}\prod_{n\geq1}\left(1-\frac{z^2}{{{j^{\prime}}_{\nu,n}^2}(q)}\right)
\end{equation*}
and
\begin{equation*}
\frac{dJ_{\nu}^{(3)}(z;q)}{dz}=\frac{\nu {z}^{{\nu}-1}}{c_{\nu}(q)}\prod_{n\geq1}\left(1-\frac{z^2}{{{l^{\prime}}_{\nu,n}^2}(q)}\right)
\end{equation*}
where $ j_{\nu,n}^{\prime}(q) $ and $ l_{\nu,n}^{\prime}(q) $ are the $ n $th positive zeros of the functions $z\mapsto {dJ_{\nu}^{(2)}(z;q)}/{dz} $ and $z\mapsto{dJ_{\nu}^{(3)}(z;q)}/{dz}$.

In addition, for the derivatives of the functions 
$z\mapsto g_{\nu}^{(2)}(z;q), z\mapsto h_{\nu}^{(2)}(z;q), z\mapsto g_{\nu}^{(3)}(z;q)$ and $z\mapsto h_{\nu}^{(3)}(z;q)$, the infinite product representations has been given, respectively, in \cite[Lemma 8., p.975]{BDM} as follow:
\begin{equation}\label{2.1}
\frac{dg_{\nu}^{(2)}(z;q)}{dz}=\prod_{n\geq1}\left(1-\frac{z^2}{\alpha_{\nu,n}^2(q)}\right),
\end{equation} 
\begin{equation}\label{2.2}
\frac{dh_{\nu}^{(2)}(z;q)}{dz}=\prod_{n\geq1}\left(1-\frac{z}{\beta_{\nu,n}^2(q)}\right),
\end{equation}
\begin{equation}\label{2.3}
\frac{dg_{\nu}^{(3)}(z;q)}{dz}=\prod_{n\geq1}\left(1-\frac{z^2}{\gamma_{\nu,n}^2(q)}\right)
\end{equation}
and
\begin{equation}\label{2.4}
\frac{dh_{\nu}^{(3)}(z;q)}{dz}=\prod_{n\geq1}\left(1-\frac{z}{\delta_{\nu,n}^2(q)}\right),
\end{equation}
where $\alpha_{\nu,n}(q)$ and $\beta_{\nu,n}(q)$ are the $n$th positive zeros of $ z\mapsto z.{dJ_{\nu}^{(2)}(z;q)}/{dz}+(1-\nu)J_{\nu}^{(2)}(z;q)$ and $ z\mapsto z.{dJ_{\nu}^{(2)}(z;q)}/{dz}+(2-\nu)J_{\nu}^{(2)}(z;q)$, while $\gamma_{\nu,n}(q)$ and $\delta_{\nu,n}(q)$ are the $n$th positive zeros of $ z\mapsto z.{dJ_{\nu}^{(3)}(z;q)}/{dz}+(1-\nu)J_{\nu}^{(3)}(z;q)$ and $ z\mapsto z.{dJ_{\nu}^{(3)}(z;q)}/{dz}+(2-\nu)J_{\nu}^{(3)}(z;q)$.

Finally, it is known from \cite[Lemma 9., p. 975]{BDM} that, between any two consecutive roots of the function $z\mapsto J_{\nu}^{(s)}(z;q)$ the function $ z\mapsto {dJ_{\nu}^{(s)}(z;q)}/{dz}$ has precisely one zero when $\nu\geq0$ and $s\in\{2,3\}$.

Now we present our first main results on the $q$-Bessel functions.

\begin{theorem}
	Let $\nu>-1, s\in\{2,3\} \text{ and } q\in(0,1).$ The following assertions are true.

	\item[\bf a.] Suppose that $ \nu>0$. Then, the radius of uniform convexity of the function $z\mapsto f_{\nu}^{(s)}(z;q)$ is the smallest positive root of the equation $$1+2r\frac{\left(f_{\nu}^{(s)}(r;q)\right)^{\prime\prime}}{\left(f_{\nu}^{(s)}(r;q)\right)^{\prime}}=0.$$
	
	\item[\bf b.] The radius of uniform convexity of the function $z\mapsto g_{\nu}^{(s)}(z;q)$ is the smallest positive root of the equation $$(2\nu-1)(\nu-1)J_{\nu}^{(s)}(r;q)+(5-4\nu)r\left(J_{\nu}^{(s)}(r;q)\right)^{\prime}+2r^2\left(J_{\nu}^{(s)}(r;q)\right)^{\prime\prime}=0.$$
	
	\item[\bf c.] The radius of uniform convexity of the function $z\mapsto h_{\nu}^{(s)}(z;q)$ is the smallest positive root of the equation $$(\nu-1)(\nu-2)J_{\nu}^{(s)}(\sqrt{r};q)+(4-2\nu)\sqrt{r}\left(J_{\nu}^{(s)}(\sqrt{r};q)\right)^{\prime}+r\left(J_{\nu}^{(s)}(\sqrt{r};q)\right)^{\prime\prime}=0.$$
\end{theorem}
\begin{proof}[\bf Proof] 
	The proofs for the cases $s=2\text{ and }s=3$ are almost the same. This is why we are going to present the proof only for the case $s=2.$ 
		\item[\bf a.]   Let $j_{\nu,n}(q) $ and $j_{\nu,n}^{\prime}(q)$ be the $ n $th positive roots of the functions 
		$z\mapsto J_{\nu}^{(2)}(z;q) $ and $z\mapsto {dJ_{\nu}^{(2)}(z;q)}/{dz} $, respectively. In \cite[p. 979]{BDM}, it was shown that the following equality is valid 
		 $$1+z\frac{\left(f_{\nu}^{(2)}(z;q)\right)^{\prime\prime}}{\left(f_{\nu}^{(2)}(z;q)\right)^{\prime}}=1-\left(\frac{1}{\nu}-1\right)\sum_{n\geq1}\frac{2z^2}{{j}^2_{\nu,n}(q)-z^2}-\sum_{n\geq1}\frac{2z^2}{{j^{\prime}}^2_{\nu,n}(q)-z^2}.$$ In the first step of our proof we consider the case $ \nu\geq1 $. We know that the zeros of Jackson's second and third $q$-Bessel functions are all real when $\nu>-1,$ according to \cite{ismail1,koelink}. Also, it is known from \cite[Lemma 9., p. 975]{BDM} that the zeros of the functions  $z\mapsto J_{\nu}^{(s)}(z;q)$ and $z\mapsto {dJ_{\nu}^{(s)}(z;q)}/{dz}$ are interlace. Here it is important to mention that the non-negative smallest zero is $z=0$ for Jackson's second and third $q$-Bessel functions. By taking $\lambda=1-\frac{1}{\nu}$ in the inequality \eqref{1.3} we have
		 $$\Re\left(\frac{2z^2}{{j^{\prime}}^2_{\nu,n}(q)-z^2}-\left(1-\frac{1}{\nu}\right)\frac{2z^2}{{j}^2_{\nu,n}(q)-z^2}\right)\leq\left(\frac{2r^2}{{j^{\prime}}^2_{\nu,n}(q)-r^2}-\left(1-\frac{1}{\nu}\right)\frac{2r^2}{{j}^2_{\nu,n}(q)-r^2}\right),$$ for $\left|z\right|\leq r<j_{\nu,1}^{\prime}(q)<j_{\nu,1}(q)$ and so, we get that
		 \begin{equation}\label{2.5}
		 \Re\left(1+z\frac{\left(f_{\nu}^{(2)}(z;q)\right)^{\prime\prime}}{\left(f_{\nu}^{(2)}(z;q)\right)^{\prime}}\right)\geq1+r\frac{\left(f_{\nu}^{(2)}(r;q)\right)^{\prime\prime}}{\left(f_{\nu}^{(2)}(r;q)\right)^{\prime}}.
		 \end{equation}
		 On the other hand, the inequality \eqref{1.2} implies that
		 $$\left|\frac{2z^2}{{j^{\prime}}^2_{\nu,n}(q)-z^2}-\left(1-\frac{1}{\nu}\right)\frac{2z^2}{{j}^2_{\nu,n}(q)-z^2}\right|\leq\frac{2r^2}{{j^{\prime}}^2_{\nu,n}(q)-r^2}-\left(1-\frac{1}{\nu}\right)\frac{2r^2}{{j}^2_{\nu,n}(q)-r^2},$$ where $\left|z\right|\leq r<j_{\nu,1}^{\prime}(q)<j_{\nu,1}(q)$. Therefore, we obtain that
		 \begin{equation}\label{2.6}
		 \left|z\frac{\left(f_{\nu}^{(2)}(z;q)\right)^{\prime\prime}}{\left(f_{\nu}^{(2)}(z;q)\right)^{\prime}}\right|\leq-r\frac{\left(f_{\nu}^{(2)}(r;q)\right)^{\prime\prime}}{\left(f_{\nu}^{(2)}(r;q)\right)^{\prime}}.
		 \end{equation}
		As second step, one can easily show that the inequalities \eqref{2.5} and \eqref{2.6} hold for $\nu\in(0,1).$ Clearly, by considering the inequality \eqref{1.4} we can write that
		$$\Re\left(\frac{2z^2}{{j^{\prime}}^2_{\nu,n}(q)-z^2}\right)\leq\left|\frac{2z^2}{{j^{\prime}}^2_{\nu,n}(q)-z^2}\right|\leq\frac{2r^2}{{j^{\prime}}^2_{\nu,n}(q)-r^2}$$
		and
		$$\Re\left(\frac{2z^2}{{j}^2_{\nu,n}(q)-z^2}\right)\leq\left|\frac{2z^2}{{j}^2_{\nu,n}(q)-z^2}\right|\leq\frac{2r^2}{{j}^2_{\nu,n}(q)-r^2}$$
		for $\left|z\right|\leq r<j_{\nu,1}^{\prime}(q)<j_{\nu,1}(q)$. Since $\frac{1}{\nu}-1>0,$ the above last two inequalities imply that the inequalities \eqref{2.5} and \eqref{2.6} hold true. Consequently, using these two inequalities yields that
		\begin{equation}\label{2.7}
		\Re\left(1+z\frac{\left(f_{\nu}^{(2)}(z;q)\right)^{\prime\prime}}{\left(f_{\nu}^{(2)}(z;q)\right)^{\prime}}\right)-\left|z\frac{\left(f_{\nu}^{(2)}(z;q)\right)^{\prime\prime}}{\left(f_{\nu}^{(2)}(z;q)\right)^{\prime}}\right|\geq1+2r\frac{\left(f_{\nu}^{(2)}(r;q)\right)^{\prime\prime}}{\left(f_{\nu}^{(2)}(r;q)\right)^{\prime}}  
		\end{equation}
		for $\left|z\right|\leq r<j_{\nu,1}^{\prime}(q).$ In \eqref{2.7}, the equality holds if and only if $z=r.$ Thus, it follws that
		$$\inf_{\left|z\right|<r}\left[\Re\left(1+z\frac{\left(f_{\nu}^{(2)}(z;q)\right)^{\prime\prime}}{\left(f_{\nu}^{(2)}(z;q)\right)^{\prime}}\right)-\left|z\frac{\left(f_{\nu}^{(2)}(z;q)\right)^{\prime\prime}}{\left(f_{\nu}^{(2)}(z;q)\right)^{\prime}}\right|\right]=1+2r\frac{\left(f_{\nu}^{(2)}(r;q)\right)^{\prime\prime}}{\left(f_{\nu}^{(2)}(r;q)\right)^{\prime}},$$ where $r\in(0,j_{\nu,1}^{\prime}(q))$. The mapping $\Phi_\nu:(0,j_{\nu,1}^{\prime}(q))\mapsto\mathbb{R}$ defined by
		$$\Phi_\nu(r)=1+2r\frac{\left(f_{\nu}^{(2)}(r;q)\right)^{\prime\prime}}{\left(f_{\nu}^{(2)}(r;q)\right)^{\prime}}=1-2\sum_{n\geq1}\left(\frac{2r^2}{{j^{\prime}}^2_{\nu,n}(q)-r^2}-\left(1-\frac{1}{\nu}\right)\frac{2r^2}{{j}^2_{\nu,n}(q)-r^2}\right)$$
		is strictly decreasing since
		$$\Phi^{\prime}_\nu(r)=-2\sum_{n\geq1}\left(\frac{4r{j^{\prime}}^2_{\nu,n}(q)}{\left({j^{\prime}}^2_{\nu,n}(q)-r^2\right)^2}-\left(1-\frac{1}{\nu}\right)\frac{4r{j}^2_{\nu,n}(q)}{\left({j}^2_{\nu,n}(q)-r^2\right)^2}\right)<0$$ for $r\in(0,j_{\nu,1}^{\prime}(q))$. Also, we have the following limits
			$$\lim_{r\searrow0}\Phi_\nu(r)=1\text{ and }\lim_{r\nearrow{j_{\nu,1}^{\prime}(q)}}\Phi_\nu(r)=-\infty.$$ As a result of this, we can say that the equation $$1+2r\frac{\left(f_{\nu}^{(2)}(r;q)\right)^{\prime\prime}}{\left(f_{\nu}^{(2)}(r;q)\right)^{\prime}}=0$$ has a unique root $ r_0 $ in the interval $(0,j_{\nu,1}^{\prime}(q))$ which is the radius of uniform convexity $r_0=r^{uc}\left(f_{\nu}^{(2)}(z;q)\right)$ of the function $z\mapsto f_{\nu}^{(2)}(z;q).$
		\item[\bf b.] By using logarithmic derivative of the function $ z\mapsto {dg_{\nu}^{(2)}(z;q)}/{dz} $ which is given by \eqref{2.1} we get that
		\begin{equation}\label{2.8}
		z\frac{\left(g_{\nu}^{(2)}(z;q)\right)^{\prime\prime}}{\left(g_{\nu}^{(2)}(z;q)\right)^{\prime}}=-\sum_{n\geq1}\frac{2z^2}{{\alpha}^2_{\nu,n}(q)-z^2}
		\end{equation}
		and
		\begin{equation}\label{2.9}
		1+z\frac{\left(g_{\nu}^{(2)}(z;q)\right)^{\prime\prime}}{\left(g_{\nu}^{(2)}(z;q)\right)^{\prime}}=1-\sum_{n\geq1}\frac{2z^2}{{\alpha}^2_{\nu,n}(q)-z^2}.
		\end{equation}
		Now, for $ \left|z\right|\leq{r}<\alpha_{\nu,1}(q) $, using the inequality \eqref{1.4} in the equalities \eqref{2.9} and \eqref{2.8}, respectively, imply that
		\begin{equation}\label{2.10}
		\Re\left(1+z\frac{\left(g_{\nu}^{(2)}(z;q)\right)^{\prime\prime}}{\left(g_{\nu}^{(2)}(z;q)\right)^{\prime}}\right)\geq1+r\frac{\left(g_{\nu}^{(2)}(r;q)\right)^{\prime\prime}}{\left(g_{\nu}^{(2)}(r;q)\right)^{\prime}}
		\end{equation}
		and
		\begin{equation}\label{2.11}
		\left|z\frac{\left(g_{\nu}^{(2)}(z;q)\right)^{\prime\prime}}{\left(g_{\nu}^{(2)}(z;q)\right)^{\prime}}\right|\leq-r\frac{\left(g_{\nu}^{(2)}(r;q)\right)^{\prime\prime}}{\left(g_{\nu}^{(2)}(r;q)\right)^{\prime}}.
		\end{equation}
		From the inequalities \eqref{2.10} and \eqref{2.11}, we deduce
		\begin{equation}\label{2.12}
		\Re\left(1+z\frac{\left(g_{\nu}^{(2)}(z;q)\right)^{\prime\prime}}{\left(g_{\nu}^{(2)}(z;q)\right)^{\prime}}\right)-\left|z\frac{\left(g_{\nu}^{(2)}(z;q)\right)^{\prime\prime}}{\left(g_{\nu}^{(2)}(z;q)\right)^{\prime}}\right|\geq1+2r\frac{\left(g_{\nu}^{(2)}(r;q)\right)^{\prime\prime}}{\left(g_{\nu}^{(2)}(r;q)\right)^{\prime}}  
		\end{equation}
		for $ \left|z\right|\leq{r}<\alpha_{\nu,1}(q) $. Equality holds in \eqref{2.12} if and only if $z=r$. As a result, we have
		$$\inf_{\left|z\right|<r}\left[\Re\left(1+z\frac{\left(g_{\nu}^{(2)}(z;q)\right)^{\prime\prime}}{\left(g_{\nu}^{(2)}(z;q)\right)^{\prime}}\right)-\left|z\frac{\left(g_{\nu}^{(2)}(z;q)\right)^{\prime\prime}}{\left(g_{\nu}^{(2)}(z;q)\right)^{\prime}}\right|\right]=1+2r\frac{\left(g_{\nu}^{(2)}(r;q)\right)^{\prime\prime}}{\left(g_{\nu}^{(2)}(r;q)\right)^{\prime}},$$
		where $r\in\left(0,\alpha_{\nu,1}(q)\right).$ Now consider the function $A_\nu:\left(0,\alpha_{\nu,1}(q)\right)\mapsto\mathbb{R}$ defined by
		$$A_\nu(r)=1+2r\frac{\left(g_{\nu}^{(2)}(r;q)\right)^{\prime\prime}}{\left(g_{\nu}^{(2)}(r;q)\right)^{\prime}}=1-\sum_{n\geq1}\frac{4r^2}{{\alpha}^2_{\nu,n}(q)-r^2}.$$ The function $ A_\nu(r) $ is strictly decreasing since 
		$${A_\nu^\prime}(r)=-\sum_{n\geq1}\frac{8r{\alpha}^2_{\nu,n}(q)}{\left({\alpha}^2_{\nu,n}(q)-r^2\right)^2}<0$$ for $r\in\left(0,\alpha_{\nu,1}(q)\right)$ and also
		$$\lim_{r\searrow0}A_\nu(r)=1\text{ and }\lim_{r\nearrow\alpha_{\nu,1}(q)}A_\nu(r)=-\infty.$$ Therefore, the equation
		\begin{equation}\label{2.13}
		1+2r\frac{\left(g_{\nu}^{(2)}(r;q)\right)^{\prime\prime}}{\left(g_{\nu}^{(2)}(r;q)\right)^{\prime}}=0
		\end{equation}
		has a unique root $r_1\in\left(0,\alpha_{\nu,1}(q)\right)$ and $r_1=r^{uc}\left(g_{\nu}^{(2)}(z;q)\right).$ By using the first and second derivatives of the function $z\mapsto g_{\nu}^{(2)}(z;q)$, one can easily see that the equation \eqref{2.13} is equivalent to $$(2\nu-1)(\nu-1)J_{\nu}^{(2)}(r;q)+(5-4\nu)r\left(J_{\nu}^{(2)}(r;q)\right)^{\prime}+2r^2\left(J_{\nu}^{(2)}(r;q)\right)^{\prime\prime}=0.$$ So, the proof is completed.
		\item[\bf c.] The proof of this part can be done similar manner. Logarithmic derivative of the function $ z\mapsto {dh_{\nu}^{(2)}(z;q)}/{dz} $ which is given by \eqref{2.2} implies that
		\begin{equation}\label{2.14}
		z\frac{\left(h_{\nu}^{(2)}(z;q)\right)^{\prime\prime}}{\left(h_{\nu}^{(2)}(z;q)\right)^{\prime}}=-\sum_{n\geq1}\frac{z}{{\beta}^2_{\nu,n}(q)-z}
		\end{equation}
		and
		\begin{equation}\label{2.15}
		1+z\frac{\left(h_{\nu}^{(2)}(z;q)\right)^{\prime\prime}}{\left(h_{\nu}^{(2)}(z;q)\right)^{\prime}}=1-\sum_{n\geq1}\frac{z}{{\beta}^2_{\nu,n}(q)-z}.
		\end{equation}
		Now, for $ \left|z\right|\leq{r}<{\beta}^2_{\nu,1}(q)$, by using the inequality \eqref{1.4} in the equalities \eqref{2.15} and \eqref{2.14}, respectively, we get that
		\begin{equation}\label{2.16}
		\Re\left(1+z\frac{\left(h_{\nu}^{(2)}(z;q)\right)^{\prime\prime}}{\left(h_{\nu}^{(2)}(z;q)\right)^{\prime}}\right)\geq1+r\frac{\left(h_{\nu}^{(2)}(r;q)\right)^{\prime\prime}}{\left(h_{\nu}^{(2)}(r;q)\right)^{\prime}}
		\end{equation}
		and
		\begin{equation}\label{2.17}
		\left|z\frac{\left(h_{\nu}^{(2)}(z;q)\right)^{\prime\prime}}{\left(h_{\nu}^{(2)}(z;q)\right)^{\prime}}\right|\leq-r\frac{\left(h_{\nu}^{(2)}(r;q)\right)^{\prime\prime}}{\left(h_{\nu}^{(2)}(r;q)\right)^{\prime}}.
		\end{equation}
		Now summarizing the inequalities \eqref{2.16} and \eqref{2.17}, we obtain
		\begin{equation}\label{2.18}
		\Re\left(1+z\frac{\left(h_{\nu}^{(2)}(z;q)\right)^{\prime\prime}}{\left(h_{\nu}^{(2)}(z;q)\right)^{\prime}}\right)-\left|z\frac{\left(h_{\nu}^{(2)}(z;q)\right)^{\prime\prime}}{\left(h_{\nu}^{(2)}(z;q)\right)^{\prime}}\right|\geq1+2r\frac{\left(h_{\nu}^{(2)}(r;q)\right)^{\prime\prime}}{\left(h_{\nu}^{(2)}(r;q)\right)^{\prime}}  
		\end{equation}
		for $ \left|z\right|\leq{r}<{\beta}^2_{\nu,1}(q)$. Equality holds in \eqref{2.18} if and only if $z=r$. Finally, we have
		$$\inf_{\left|z\right|<r}\left[\Re\left(1+z\frac{\left(h_{\nu}^{(2)}(z;q)\right)^{\prime\prime}}{\left(h_{\nu}^{(2)}(z;q)\right)^{\prime}}\right)-\left|z\frac{\left(h_{\nu}^{(2)}(z;q)\right)^{\prime\prime}}{\left(h_{\nu}^{(2)}(z;q)\right)^{\prime}}\right|\right]=1+2r\frac{\left(h_{\nu}^{(2)}(r;q)\right)^{\prime\prime}}{\left(h_{\nu}^{(2)}(r;q)\right)^{\prime}},$$
		where $r\in\left(0,{\beta}^2_{\nu,1}(q)\right).$ Now consider the function $B_\nu:\left(0,{\beta}^2_{\nu,1}(q)\right)\mapsto\mathbb{R}$ defined by
		$$B_\nu(r)=1+2r\frac{\left(h_{\nu}^{(2)}(r;q)\right)^{\prime\prime}}{\left(h_{\nu}^{(2)}(r;q)\right)^{\prime}}=1-\sum_{n\geq1}\frac{2r}{{\beta}^2_{\nu,n}(q)-r}.$$ The function $ B_\nu(r) $ is strictly decreasing since 
		$${B_\nu^\prime}(r)=-\sum_{n\geq1}\frac{2{\beta}^2_{\nu,n}(q)}{\left({\beta}^2_{\nu,n}(q)-r\right)^2}<0$$ for $r\in\left(0,{\beta}^2_{\nu,1}(q)\right)$ and also
		$$\lim_{r\searrow0}B_\nu(r)=1\text{ and }\lim_{r\nearrow{\beta}^2_{\nu,1}(q)}B_\nu(r)=-\infty.$$ As a result, the equation
		\begin{equation}\label{2.19}
		1+2r\frac{\left(h_{\nu}^{(2)}(r;q)\right)^{\prime\prime}}{\left(h_{\nu}^{(2)}(r;q)\right)^{\prime}}=0
		\end{equation}
		has a unique root $r_2\in\left(0,{\beta}^2_{\nu,1}(q)\right)$ and $r_2=r^{uc}\left(h_{\nu}^{(2)}(z;q)\right).$ By considering the first and second derivatives of the function $z\mapsto h_{\nu}^{(2)}(z;q)$, we can easily obtain that the equation \eqref{2.19} is equivalent to $$(\nu-1)(\nu-2)J_{\nu}^{(2)}(\sqrt{r};q)+(4-2\nu)\sqrt{r}\left(J_{\nu}^{(2)}(\sqrt{r};q)\right)^{\prime}+r\left(J_{\nu}^{(2)}(\sqrt{r};q)\right)^{\prime\prime}=0,$$ which is desired.
\end{proof}

\subsection{Uniform Convexity of some normalized Wright Functions}
In this subsection, we will focus on the function
$$\phi(\rho,\beta,z)=\sum_{n\geq0}\frac{z^n}{n!\Gamma(n\rho+\beta)}  \text{ \ \ \ \ } (\rho>-1 \text{ \ \and\ \ }z,\beta\in\mathbb{C}) $$
named after the British mathematician E.M. Wright. It is well known that this function was introduced by him for the first time in the case $\rho>0$ in connection with his investigations on the asymptotic theory of partitions \cite{wright}. 

\bigskip
From \cite[Lemma 1]{btk} we know that under the conditions $\rho>0$ and $\beta>0,$ the function $z\mapsto\lambda_{\rho,\beta}(z)=\phi(\rho,\beta,-z^2)$ has infinitely many zeros which are all real. Thus, due to the Hadamard factorization theorem, the expression $\lambda_{\rho,\beta}(z)$ can be written as
$$\Gamma(\beta)\lambda_{\rho,\beta}(z)=\prod_{n\geq1}\left(1-\frac{z^2}{\lambda^2_{\rho,\beta,n}}\right)$$
where $\lambda_{\rho,\beta,n}$ stands for the $n$th positive zero of the function $\lambda_{\rho,\beta}(z)$ (or the positive real zeros of the function $\Psi_{\rho,\beta}$). Moreover, let $\zeta_{\rho,\beta,n}'$ denote the $n$th positive zero of $\Psi_{\rho,\beta}',$ where $\Psi_{\rho,\beta}(z)=z^{\beta}\lambda_{\rho,\beta}(z),$ then the zeros satisfy the chain of inequalities
$$\zeta^{\prime}_{\rho,\beta,1}<\zeta_{\rho,\beta,1}<\zeta^{\prime}_{\rho,\beta,2}<\zeta_{\rho,\beta,2}<{\dots} .$$

One can easily see that the function $z\mapsto\phi(\rho,\beta,-z^2)$ do not belong to $\mathcal{A}$, and thus first we perform some natural normalization. We define three functions originating $\phi(\rho,\beta,.)$: 
\begin{equation*}
f_{\rho,\beta}(z)=\left(z^{\beta}\Gamma(\beta)\phi(\rho,\beta,-z^2)\right)^{\frac{1}{\beta}},
\end{equation*}
\begin{equation*}
g_{\rho,\beta}(z)=z\Gamma(\beta)\phi(\rho,\beta,-z^2),
\end{equation*}
\begin{equation*}
h_{\rho,\beta}(z)=z\Gamma(\beta)\phi(\rho,\beta,-z).
\end{equation*}
Clearly these functions are contained in the class $\mathcal{A}$.

 Now, we would like to present our results regarding the uniform convexity of the functions $f_{\rho,\beta}$, $g_{\rho,\beta}$ and $h_{\rho,\beta}.$
 
\begin{theorem}
	Let $\rho>0$ and $\beta>0$.
	\item{\bf a.} The radius of uniform convexity of the function $f_{\rho,\beta}$ is the smallest positive root of the equation
	$$1+2r\frac{\Psi^{\prime \prime}_{\rho,\beta}(r)}{\Psi^{\prime}_{\rho,\beta}(r)}+2\left( \frac{1}{\beta}-1\right) \frac{r\Psi^{\prime}_{\rho,\beta}(r)}{\Psi_{\rho,\beta}(r)}=0,$$
	where $\Psi_{\rho,\beta}(z)=z^{\beta}\lambda_{\rho,\beta}(z).$
	\item{\bf b.} The radius of uniform convexity of the function $g_{\rho,\beta}$ is the smallest positive root of the equation $$1+2r\frac{g_{\rho,\beta}''(r)}{g_{\rho,\beta}'(r)}=0.$$
	\item{\bf c.} The radius of uniform convexity of the function $h_{\rho,\beta}$ is the smallest positive root of the equation $$1+2r\frac{h_{\rho,\beta}''(r)}{h_{\rho,\beta}'(r)}=0.$$
\end{theorem}
\begin{proof}
	{\bf a.} Let $\zeta_{\rho,\beta,n}$ and $\zeta^{\prime}_{\rho,\beta,n}$ be the $n$th positive roots of $\Psi_{\rho,\beta}$ and $\Psi^{\prime}_{\rho,\beta},$ respectively. In \cite[Theorem 5]{btk} the following equality was demonstrated:
	$$1+\frac{zf^{\prime \prime}_{\rho,\beta}(z)}{f^{\prime }_{\rho,\beta}(z)}=1-\left( \frac{1}{\beta}-1\right) \sum_{n\geq1}\frac{2z^{2}}{\zeta^{2}_{\rho,\beta,n}-z^{2}}-\sum_{n\geq1}\frac{2z^{2}}{\zeta^{\prime 2}_{\rho,\beta,n}-z^{2}}.$$
	In order to prove the theorem we need to investigate two different cases such as $\beta \in (0,1]$ and $\beta>1$. First suppose $\beta \in (0,1]$. In this case, with the help of \eqref{1.4} for  $\beta\in (0,1]$,  we deduce that the inequality
	\begin{align}
	\Re\left(1+\frac{zf^{\prime \prime}_{\rho,\beta}(z)}{f^{\prime }_{\rho,\beta}(z)}\right)&\geq 1-\left( \frac{1}{\beta}-1\right) \sum_{n\geq1}\frac{2r^{2}}{\zeta^{2}_{\rho,\beta,n}-r^{2}}-\sum_{n\geq1}\frac{2r^{2}}{\zeta^{\prime 2}_{\rho,\beta,n}-r^{2}}\label{2.20}\\&=1+\frac{rf^{\prime \prime}_{\rho,\beta}(r)}{f^{\prime }_{\rho,\beta}(r)}, \quad |z|\leq r<\zeta^{\prime}_{\rho,\beta,1}<\zeta_{\rho,\beta,1}. \nonumber
	\end{align}
	holds true for $\left| z\right|=r$. Moreover, in view of \eqref{1.4}, we get
	\begin{align}
	\left|\frac{zf^{\prime \prime}_{\rho,\beta}(z)}{f^{\prime }_{\rho,\beta}(z)} \right|&=\left| \sum_{n\geq1}\frac{2z^{2}}{\zeta^{\prime 2}_{\rho,\beta,n}-z^{2}}+\left( \frac{1}{\beta}-1\right)\sum_{n\geq1}\frac{2z^{2}}{\zeta^{2}_{\rho,\beta,n}-z^{2}}\right| \label{2.21}\\
	&\leq \sum_{n\geq1}\left| \left(\frac{2z^{2}}{\zeta^{\prime 2}_{\rho,\beta,n}-z^{2}}+\left( \frac{1}{\beta}-1\right)\frac{2z^{2}}{\zeta^{2}_{\rho,\beta,n}-z^{2}} \right)\right| \nonumber\\
	&\leq \sum_{n\geq1}\left(\frac{2r^{2}}{\zeta^{\prime 2}_{\rho,\beta,n}-r^{2}}+\left( \frac{1}{\beta}-1\right)\frac{2r^{2}}{\zeta^{2}_{\rho,\beta,n}-r^{2}} \right) \nonumber\\
	&=-\frac{rf^{\prime \prime}_{\rho,\beta}(r)}{f^{\prime }_{\rho,\beta}(r)} \nonumber
	\end{align}
	where $|z|\leq r<\zeta^{\prime}_{\rho,\beta,1}<\zeta_{\rho,\beta,1}.$ On the other hand, in view of the inequality \eqref{1.3} we obtain that \eqref{2.20} and \eqref{2.21} are also valid when $\beta\geq1$ for all  $z\in (0,{\zeta^{\prime}_{\rho,\beta,1}}).$ 
	Here we used tacitly that the zeros of $\zeta_{\rho,\beta,n}$ and $\zeta^{\prime}_{\rho,\beta,n}$ interlace as aforesaid , that is, we have $\zeta^{\prime}_{\rho,\beta,1}<\zeta_{\rho,\beta,1}.$
	Eventually, thanks to \eqref{2.20} and \eqref{2.21} we arrive at
	\begin{equation}\label{2.22}
	\Re\left(1+\frac{zf^{\prime \prime}_{\rho,\beta}(z)}{f^{\prime }_{\rho,\beta}(z)}\right)-\left|\frac{zf^{\prime\prime}_{\rho,\beta}(z)}{f^{\prime }_{\rho,\beta}(z)} \right|\geq1+2r\frac{f^{\prime \prime}_{\rho,\beta}(r)}{f^{\prime }_{\rho,\beta}(r)}, \quad \left| z\right| \leq r<\zeta^{\prime}_{\rho,\beta,1}.
	\end{equation}
	Due to the minimum principle for harmonic functions, equality holds if and only if $z=r$. Now, the above deduced inequalities imply for $r\in(0,\zeta^{\prime}_{\rho,\beta,1})$
	$$\inf_{\left| z\right|< r}\left\lbrace \Re\left(1+\frac{zf^{\prime \prime}_{\rho,\beta}(z)}{f^{\prime }_{\rho,\beta}(z)}\right)-\left|\frac{zf^{\prime \prime}_{\rho,\beta}(z)}{f^{\prime }_{\rho,\beta}(z)} \right|\right\rbrace =1+2r\frac{f^{\prime \prime}_{\rho,\beta}(r)}{f^{\prime }_{\rho,\beta}(r)}.$$
	On the other hand, the function $u_{\rho,\beta}:(0,\zeta^{\prime}_{\rho,\beta,1})\rightarrow\mathbb{R},$ defined by
	$$u_{\rho,\beta}(r)=1+2\frac{rf^{\prime\prime}_{\rho,\beta}(r)}{f^{\prime}_{\rho,\beta}(r)}=1-2\sum_{n\geq1}\left(\frac{2r^{2}}{\zeta^{\prime2}_{\rho,\beta,n}-r^{2}}-\left(1-\frac{1}{\beta}\right)\frac{2r^{2}}{\zeta^{2}_{\rho,\beta,n}-r^{2}} \right),$$
	is strictly decreasing when $\beta\in(0,1].$ Moreover, it is also strictly decreasing when $\beta>1$ since
	\begin{align*}
	u^{\prime}_{\rho,\beta}(r)&=-\left(\frac{1}{\beta}-1\right)\sum_{n\geq1}\frac{8r\zeta^{2}_{\rho,\beta,n}}{(\zeta^{2}_{\rho,\beta,n}-r^{2})^{2}}-\sum_{n\geq1}\frac{8r\zeta^{\prime2}_{\rho,\beta,n}}{(\zeta^{\prime 2}_{\rho,\beta,n}-r^{2})^{2}}
	\\&<\sum_{n\geq1}\frac{8r\zeta^{2}_{\rho,\beta,n}}{(\zeta^{2}_{\rho,\beta,n}-r^{2})^{2}}-\sum_{n\geq1}\frac{8r\zeta^{\prime 2}_{\rho,\beta,n}}{(\zeta^{\prime 2}_{\rho,\beta,n}-r^{2})^{2}}<0
	\end{align*}
	for $r\in(0,\zeta^{\prime}_{\rho,\beta,1}).$ Observe also that $$\lim_{r\searrow0}u_{\rho,\beta}(r)=1 \text{ and }\lim_{r\nearrow\zeta^{\prime}_{\rho,\beta,1}}u_{\rho,\beta}(r)=-\infty.$$ Thus it follows that the equation $$1+2r\frac{rf^{\prime \prime}_{\rho,\beta}(r)}{f^{\prime }_{\rho,\beta}(r)}=0$$ has a unique root $r_{3}\in(0,\zeta^{\prime}_{\rho,\beta,1})$ and $r_{3}=r^{uc}(f_{\rho,\beta})$.
	
	{\bf b.} Let $\vartheta_{\rho,\beta,n}$ be the $n$th positive zero of the function $g^{\prime}_{\rho,\beta}(z)$. In \cite[Theorem 5]{btk} the following equality was proven:
	\begin{equation}\label{2.23}
	1+\frac{zg''_{\rho,\beta}(z)}{g'_{\rho,\beta}(z)}=1-\sum_{n\geq1}\frac{2z^2}{\vartheta_{\rho,\beta,n}^2-z^2}.
	\end{equation}
	As a result of this equality, the inequality
	\begin{equation}\label{2.24}
	\Re\left(1+\frac{zg''_{\rho,\beta}(z)}{g'_{\rho,\beta}(z)}\right)\geq1-\sum_{n\geq1}\frac{2r^2}{\vartheta_{\rho,\beta,n}^2-r^2},\quad \left| z\right|\leq r<\vartheta_{\rho,\beta,1} 
	\end{equation}
	was shown in \cite{btk}. From equality \eqref{2.23} we arrive at
	\begin{align}
	\left|\frac{zg''_{\rho,\beta}(z)}{g'_{\rho,\beta}(z)} \right|&=\left|\sum_{n\geq1}\frac{2z^2}{\vartheta_{\rho,\beta,n}^2-z^2} \right|\leq \sum_{n\geq1}\left| \frac{2z^2}{\vartheta_{\rho,\beta,n}^2-z^2}\right| \leq \sum_{n\geq1}\frac{2r^2}{\vartheta_{\rho,\beta,n}^2-r^2}  \label{2.25}\\
	&=-\frac{rg''_{\rho,\beta}(r)}{g'_{\rho,\beta}(r)}, \quad \left| z\right|\leq r<\vartheta_{\rho,\beta,1} \nonumber. 
	\end{align}
	By using the inequalities \eqref{2.24} and \eqref{2.25} we obtain
	\begin{equation*}
	\Re\left(1+\frac{zg_{\rho,\beta}''(z)}{g_{\rho,\beta}'(z)}\right) -\left|\frac{zg''_{\rho,\beta}(z)}{g'_{\rho,\beta}(z)} \right|\geq 1+2r\frac{g_{\rho,\beta}''(r)}{g_{\rho,\beta}'(r)} \quad \left| z\right|<r<\vartheta_{\rho,\beta,1} .
	\end{equation*}
	Owing to the minimum principle for harmonic functions, equality holds if and only if $z=r$. Thus, for $r\in(0,\vartheta_{\rho,\beta,1})$ we get
	$$\inf_{\left|z\right|<r}\left\{\Re\left(1+\frac{zg_{\rho,\beta}''(z)}{g_{\rho,\beta}'(z)}\right)-\left|\frac{zg''_{\rho,\beta}(z)}{g'_{\rho,\beta}(z)}\right|\right\}=1+2r\frac{g_{\rho,\beta}''(r)}{g_{\rho,\beta}'(r)}.$$
	The function
	$v_{\rho,\beta}:(0,\vartheta_{\rho,\beta,1})\to\mathbb{R},$ defined by
	$$v_{\rho,\beta}(r)=1+2r\frac{g_{\rho,\beta}''(r)}{g_{\rho,\beta}'(r)},$$ is strictly decreasing
	and $$\lim_{r\searrow0}v_{\rho,\beta}(r)=1, \ \ \ \lim_{r\nearrow\vartheta_{\rho,\beta,1}}v_{\rho,\beta}(r)=-\infty.$$
	Consequently, the equation
	$$1+2r\frac{g_{\rho,\beta}''(r)}{g_{\rho,\beta}'(r)}=0$$ has a unique root $r_4$ in $(0,\vartheta_{\rho,\beta,1})$, and $r_{4}=r^{uc}(g_{\rho,\beta})$.
	
	{\bf c.} Let $\tau_{\rho,\beta,n}$ denote the $n$th positive zero of the function $h_{\rho,\beta}(z)$. In \cite[Theorem 5]{btk} the following equation was obtained
	\begin{equation}\label{2.26}
	\frac{zh_{\rho,\beta}''(z)}{h_{\rho,\beta}'(z)}=-\sum_{n\geq1}\frac{z}{\tau_{\rho,\beta,n}-z},
	\end{equation}
	and, in the same paper, with the help of \eqref{2.26} the following inequality was given
	\begin{equation}\label{2.27}
	\Re\left(1+\frac{zh_{\rho,\beta}''(z)}{h_{\rho,\beta}'(z)}\right)\geq 1+ \frac{rh_{\rho,\beta}''(r)}{h_{\rho,\beta}'(r)}, \quad \left|z\right|<r<\tau_{\rho,\beta,1}<\lambda_{\rho,\beta,1}.
	\end{equation}
	From \eqref{2.26} we get
	\begin{align}
	\left|\frac{zh_{\rho,\beta}''(z)}{h_{\rho,\beta}'(z)}\right|&=\left|\sum_{n\geq1}\frac{z}{\tau_{\rho,\beta,n}-z} \right|\leq\sum_{n\geq1}\left| \frac{z}{\tau_{\rho,\beta,n}-z} \right|\leq\sum_{n\geq1}\frac{r}{\tau_{\rho,\beta,n}-r}\label{2.28} \\&=-\frac{rh_{\rho,\beta}''(r)}{h_{\rho,\beta}'(r)},\quad \left| z\right|<r<\tau_{\rho,\beta,1}. \nonumber 
	\end{align}
	From inequality \eqref{2.27} and \eqref{2.28} we deduce
	\begin{equation*}
	\Re\left(1+\frac{zh_{\rho,\beta}''(z)}{h_{\rho,\beta}'(z)}\right)-\left| \frac{zh_{\rho,\beta}''(z)}{h_{\rho,\beta}'(z)}\right|\geq 1+2\frac{rh_{\rho,\beta}''(r)}{h_{\rho,\beta}'(r)} \quad \left| z\right|<r<\tau_{\rho,\beta,1}.
	\end{equation*}
	Due to the minimum principle for harmonic functions, equality holds if and only if $z=r$. Thus we obtain
	$$\inf_{\left|z\right|<r}\left\{\Re\left(1+\frac{zh_{\rho,\beta}''(z)}{h_{\rho,\beta}'(z)}\right)-\left|\frac{zh_{\rho,\beta}''(z)}{h_{\rho,\beta}'(z)}\right|\right\}=1+2r\frac{h_{\rho,\beta}''(r)}{h_{\rho,\beta}'(r)},$$
	for every $r\in (0,\tau_{\rho,\beta,1})$. Since the function $w_{\rho,\beta}(r):(0,\tau_{\rho,\beta,1})\rightarrow \mathbb{R}$ defined by $$w_{\rho,\beta}(r)=1+2r\frac{h_{\rho,\beta}''(r)}{h_{\rho,\beta}'(r)}=1-\sum_{n\geq1}\frac{2r}{\tau_{\rho,\beta,n}-r}$$
	is strictly decreasing on $(0,\tau_{\rho,\beta,1})$, and
	$$\lim_{r\searrow0}w_{\rho,\beta}(r)=1, \   \   \  \lim_{r\nearrow\tau_{\rho,\beta,1}}w_{\rho,\beta}(r)=-\infty.$$
	It follows that the equation $w_{\rho,\beta}(r)=0$ has a unique root $r_{5}\in (0,\tau_{\rho,\beta,1})$, and this root is the radius of uniform convexity.
\end{proof}

\end{document}